\documentclass[11pt]{article}

\input{amssym.def}
\input{amssym}
\usepackage{eufrak}

\newtheorem{theorem}{Theorem}

\newtheorem{conjecture}[theorem]{Conjecture}
\newtheorem{corollary}[theorem]{Corollary}
\newtheorem{definition}[theorem]{Definition}
\newtheorem{example}[theorem]{Example}

\newtheorem{proposition}[theorem]{Proposition}
\newtheorem{remark}[theorem]{Remark}

\def\pa{{\partial}}
\def\hl{\hat{l}}
\def\hL{\hat{L}}

\def\un{\underline}
\def\qq{q^{-1}}

\def\Tr{\mathrm{Tr}}
\def\Trr{\Tr_R}

\def\hLL{\hat{{\cal L}}}
\def\LL{{\cal{L}}}
\def\MM{{\cal{M}}}
\def\CC{{\cal C}}

\def\de{\delta}
\def\De{\Delta}
\def\ot{\otimes}
\def\C{{\Bbb C}}
\def\Sym{{\rm Sym\, }}
\def\vv{V^{\otimes 2}}

\def\ov{\overline}
\def\hmu{\hat{\mu}}
\def\la{{\lambda}}

\def\be{\begin{equation}}
\def\ee{\end{equation}}

\begin{document}

\makeatletter
\renewcommand{\theequation}{{\thesection}.{\arabic{equation}}}
\@addtoreset{equation}{section} \makeatother

\title{$q$-Casimir and $q$-cut-and-join operators related\\
 to Reflection Equation Algebras}
\author{\rule{0pt}{7mm} Dimitri
Gurevich\thanks{gurevich@ihes.fr}\\
{\small\it Universit\'e Polytechnique Hauts-de-France}\\
{\small\it F-59313 Valenciennes, France}\\
{\small \it and}\\
{\small \it Interdisciplinary Scientific Center J.-V.Poncelet}\\
{\small\it Moscow 119002, Russian Federation}\\
\rule{0pt}{7mm} Varvara Petrova\\
{\small\it
National Research University Higher School of Economics,}\\
{\small\it 20 Myasnitskaya Ulitsa, Moscow 101000, Russian Federation}\\
\rule{0pt}{7mm} Pavel Saponov\thanks{Pavel.Saponov@ihep.ru}\\
{\small\it
National Research University Higher School of Economics,}\\
{\small\it 20 Myasnitskaya Ulitsa, Moscow 101000, Russian Federation}\\
{\small \it and}\\
{\small \it
Institute for High Energy Physics, NRC "Kurchatov Institute"}\\
{\small \it Protvino 142281, Russian Federation}}

\maketitle

\begin{abstract}

In this paper we are dealing with the Reflection Equation algebra $\MM(R)$, associated with a $GL_N$ type Hecke symmetry $R$. In this algebra we define the $q$-analogs of
the partial derivatives $\pa_j^i$ in generators $m_i^j$ of $\MM(R)$. The product $\hL = MD$ of two matrices $M=\|m_i^j\|$ and $D=\|\partial_i^j\|$ turns out to be a generating
matrix of a modified Reflection Equation algebra $\hLL(R)$ which is similar to the universal enveloping algebra $U(gl_N)$ in many aspects. Central elements of the modified
Reflection Equation algebra give rise to $q$-Casimir operators in a representation of $\hLL(R)$ in the algebra $\MM(R)$.

We perform a spectral analysis of the first $q$-Casimir operator and formulate a conjecture about the spectrum of the higher ones. At last, we define the normal ordering for the
$q$-differential operators and inroduce the $q$-cut-and-join operators. In several explicit examples we express some of $q$-cut-and-join operators via the $q$-Casimir
ones by analogy with the classical case.
\end{abstract}

{\bf AMS Mathematics Subject Classification, 2010:} 81R50


\section{Introduction}

Let $N^2$ elements  $m_i^j$, $1\leq i,j\leq N$, be generators of the coordinate ring of the commutative algebra $\Sym (gl_N)$. We denote by $\partial_j^i = \partial/\partial m_i^j$ the
partial derivative in the indeterminate $m_i^j$.

Consider the matrix $\hL = MD$, where $N\times N$ matrices $M = \|m_i^j\|$ and $D =\|\partial_i^j\|$ are composed of the generators and the partial derivatives respectively. Hereafter
the lower index enumerates the rows of a matrix, while the upper one --- the columns. Note, that the entries $\hl_i^j$ of the matrix $\hL$ meet the commutation relations:
$$
[\hl_i^j, \hl_k^r]=\de_k^j\hl_i^r-\de_i^r\hl_k^j
$$
and consequently generate the universal enveloping algebra $U(gl_N)$. As is well known, the elements $\Tr \hL^k$, $1\le k\le N$, generate the center $Z(U(gl_N))$ of the
enveloping algebra.

Let $\la=(\la_1,\dots ,\la_N)$ be a partition and $V_\la$ be the corresponding irreducible finite dimensional $U(gl_N)$-module. The Schur lemma implies that the central elements of
$U(gl_N)$ are represented by scalar operators on the spaces $V_\la$ (we call them the Casimir operators). The eigenvalues related to $\mathrm{Tr}\hL^k$ were computed in \cite{PP}.

Lately, the following operators attracted a considerable interest:
\be
W^\De = :\!\Tr\,\hL^{\De_1}\dots \Tr\,\hL^{\De_k}\!:
\label{cut}
\ee
where $\De=(\De_1,\dots,\De_k)$ is a partition of a positive integer and colons $:\,:$ denote the normal ordering. These operators play an important role in the Hurwitz-Kontsevich
theory (see \cite{MMN}) and are called {\it the cut-and-join operators}\footnote{Usually, these operators have some normalizing factors, motivated by representation theory of symmetric groups. We disregard these factors.}. Often, the operators $W^\De$ are expressed in terms of invariants $p_k = \Tr M^k$ and partial derivatives in
$p_k$. Namely, in this form the simplest cut-and-join operator was introduced in \cite{Go}.

The main objective  of the present note is to introduce $q$-analogs of the Casimir and cut-and-join operators for the Reflection Equation (RE) algebra related to a $GL_N$ type
Hecke symmetry $R$ and to perform the spectral analysis of some of these operators. Note that in this case the properties of the RE algebra
are very similar to those of the enveloping algebra $U(gl_N)$.

By definition, {\it a Hecke symmetry} is a linear operator $R:\,\vv\to \vv$, $\dim_{\,\C}V = N$, which satisfies the so-called braid relation:
$$
(R\ot I)(I\ot R)(R\ot I)=(I\ot R)(R\ot I)(I\ot R)
$$
and the {\it Hecke condition}:
\be
(q I-R)(\qq I+R)=0,\quad q\in \C\setminus \{\pm 1,0\}.
\label{H-cond}
\ee
Hereafter, $I$ stands for the identity operator or the unit matrix.

The well known example is the Hecke symmetry $R$ coming from the Quantum Group (QG) $U_q(sl_N)$. At the limit $q\to 1$ it turns into the usual flip $P$
(the permutation operator).

Below, we impose two additional requirements on a Hecke symmetry $R$: it should be {\it skew-invertible} and {\it even}. The first requirement enables us to define the so-called
$R$-trace $\Tr_R M$ of any $N\times N$ matrix $M$. The second requirement means that the Poincar\'e-Hilbert series, corresponding to the $R$-skew-symmetric algebra
$\Lambda_R(V)$ of the basic space $V$, is a monic polynomial. The degree $m$ of this polynomial is called {\em the rank} of $R$ and in general $m\le N={\dim}_{\,\C} V$.

If $R$ satisfies the above requirenments it is called {\it a $GL_m$ type} Hecke symmetry (see \cite{GPS1} for more detail). In the present paper we mostly deal with $m=N$,
this case includes in particular the Drinfeld-Jimbo and Cremmer-Gervais $R$-matrices.

\begin{remark}\rm
There exist a lot of examples of Hecke symmetries with the rank $m<N$ which are not deformations either of the usual flips or of the super-flips (see \cite{G}).
\end{remark}

With any $GL_N$ type Hecke symmetry $R$ we associate two versions of the RE algebras --- modified and non-modified ones. The modified RE algebra $\hLL(R)$ is a unital
associative algebra generated by entries of a matrix $\hL=\|\hl_i^j\|_{1\leq i,j \leq N}$ subject to the following system of quadratic-linear relations
\be
R\hL_1R \hL_1- \hL_1R \hL_1 R = R \hL_1- \hL_1 R,\quad \hL_1=\hL\ot I.
\label{mRE}
\ee
The non-modified RE algebra $\MM(R)$ is defined by a homogeneous quadratic system of relations on entries of the generating matrix $M=\|m_i^j\|$:
$$
RM_1R M_1- M_1R M_1 R = 0,\quad M_1 = M\otimes I.
$$
Note that for any {\it Hecke symmetry} $R$ the RE algebra $\MM(R)$ is isomorphic to the modified RE algebra $\hLL(R)$. Moreover, for a Hecke symmetry $R$ the algebras
$\MM(R)$ and $\hLL(R)$ can be treated as the same RE algebra written in two different sets of generators $m_i^j$ and $\hl_i^j$ connected with each other by a simple linear
shift:
\be
M=I-(q-\qq)\hL.
\label{iso}
\ee
Nevertheless, we shall distinguish these two versions of the RE algebra and use the different notation $\MM(R)$ and $\hLL(R)$ for them. One of the reasons is that the
substitution (\ref{iso}) is singular at the limit $q \to 1$ and, consequently, we get different classical limits $q\to 1$ for $\MM(R)$ and $\hLL(R)$. Indeed, if $R\to P$ as $q\to 1$,
the algebra $\MM(R)$ tends to $\Sym(gl_N)$, whereas the algebra $\hLL(R)$ tends to $U(gl_N)$.

For any $GL_N$ type Hecke symmetry $R$, the algebra $\hLL(R)$ (respectively $\MM(R)$) has a large center  similar to the center of the algebra $U(gl_N)$. In particular,
the center is generated by the elements $\Trr \hL^k$ (respectively $\Trr M^k$), called the power sums. Besides, it is possible to define $q$-analogs of the partial derivatives
$\pa_i^j$ in the algebra $\MM(R)$. In the full analogy with the classical case, the matrix  $\hL=MD,$ where $D=\|\pa_i^j\|$, generates the modified RE algebra
$\hLL(R)$, which is similar to the enveloping algebra\footnote{If a Hecke symmetry $R$ is a deformation of the usual flip $P$ (for instance, it comes from $U_q(sl_N)$ or
it is a Cremmer-Gervais symmetry), the algebra $\hLL(R)$ can be treated as a two parameter deformation of the algebra $\Sym(gl_N)$. The quasi-classical counterpart of this
deformation  is a Poisson pencil, generated by two brackets (see \cite{GPS1}). By contrast, the QG $U_q(gl_N)$ is isomorphic to $U(sl_N)$ as an algebra, only its coalgebraic
structure is deformed.} $U(gl_N)$.

In \cite{GPS1} we constructed $\hLL(R)$-modules which are $q$-analogs of the $U(gl(N))$-modules $V_\la$. Conjecturally, the modules $V_\lambda$ are irreducible\footnote{Note
that the RE algebras admit representations different from those defined in \cite{GPS1}. For instance, some examples of reducible but indecomposable representations were
constructed in \cite{S}.} for a generic $q$. Up to this conjecture, the central elements of the algebra $\hLL(R)$ are represented by scalar operators in these modules.

It is also true for the Casimir operators of the modified RE algebra $\hLL(R)$ represented as operators acting onto another (non-modified) copy of the RE algebra $\MM(R)$
provided the action is a $q$-analog of the $U(gl_N)$ representation by the right-invariant vector fields on the algebra $\Sym(gl_N)$. We perform the spectral analysis of the first
Casimir operator. In this analysis an important role is played by some idempotents belonging to the Hecke algebras. Using these idempotents, we define the $q$-Schur polynomials
$s_\la(M)\in \MM(R)$ and prove that they are eigenvectors of the $q$-cut-and-join operator.

This note is organized as follows. In the next section, we introduce a quantum double, giving rise to the quantum partial derivatives and to the relation $\hL=MD$ similar to the
classical one. Also, we introduce the $q$-Schur polynomials and power sums. Section 3 contains a spectral analysis of the first $q$-Casimir operator (which
coincides with the first $q$-cut-and-join operator). In particular, we show that the $q$-Schur polynomials are its eigenvectors and find the corresponding eigenvalues. In the last
section we introduce the ``eigenvalues" $\mu_i$ and $\hat \mu_i$ of the quantum matrices $M$ and $\hL$ respectively and formulate a conjecture for the spectrum of the higher
$q$-Casimir and $q$-cut-and-join operators in terms of these eigenvalues. A few examples of the $q$-cut-and-join operators and their relations with the $q$-Casimir operators are
exhibited.

The ground field is assumed to be $\C$. The parameter $q$ is assumed to be generic. In particular, this means that $q^k\not=1$ for any integer $k>1$.

\section{Quantum Doubles and Casimir operators}

Given two associative unital algebras $A$ and $B$, we assume that there exists a {\it permutation map}
$$
\sigma: A\ot B\to B\ot A
$$
preserving the algebraic structures of $A$ and $B$. If these algebras are defined via relations on their generators, this property means that the relations are preserved by $\sigma$
(for a more detailed definition see \cite{GS2}). Using the map $\sigma$, one can convert the linear space $B\otimes A$ into the unital associative algebra with the product
defined by the rule:
$$
(b_1\otimes a_1)*(b_2\otimes a_2) = (m_B\otimes m_A)\circ (b_1\otimes \sigma(a_1\otimes b_2)\otimes a_2),
$$
where $m_B$ and $m_A$ are the product operations in the algebras $B$ and $A$ respectively. Such an algebra $B\otimes A$ is called a {\it quantum double} (QD) of the algebras
$A$ and $B$ if the map $\sigma$ is different from a (super-)flip. We denote the corresponding QD as $(A,B,\sigma)$ or simply $(A,B)$.

The algebras $A\simeq 1_B\otimes A$ and $B\simeq B\otimes 1_A$ are subalgebras of the QD $(A,B)$ with respect to the above $*$ product since by definition
$\sigma(1_A\otimes b) = b\otimes 1_A$ and $\sigma(a\otimes 1_B) = 1_B\otimes a$ for any $a\in A$, $b\in B$ (see \cite{GS2}). The QD products $*$ of elements of the
subalgebra $A$ and elements of the subalgebra $B$ will be called the {\it permutation relations} induced by the map $\sigma$:
$$
(1_B\otimes a)*(b\otimes 1_A) = \sigma(a\otimes b).
$$

If the algebra $A$ is endowed with a counit $\varepsilon_A:A\to \C$ (a one-dimensional representation), it is possible to define an action $\triangleright : A\ot B\rightarrow B$ of the
algebra $A$ onto $B$ by setting
\be
a\triangleright b:=(\mathrm{Id}\ot \varepsilon_A)\sigma (a\ot b), \quad \forall\, a\in A, \, b\in B.
\label{AB-act}
\ee

A typical example of a QD is the so-called Heisenberg double \cite{IP}, in which the roles of the algebras $A$ and $B$ are played respectively by an RE algebra $\LL(R)$ and
the so-called RTT algebra, defined via the same Hecke symmetry $R$.

In \cite{GPS2, GPS3} we have introduced a QD $(A,B,\sigma)$, where $A=\hLL(R)$, $B=\MM(R)$ with the generating matrices $\hL$ and $M$ respectively and the permutation map
\be
\sigma:\,R\hL_1 R M_1\to M_1R \hL_1 R^{-1}+RM_11_A. \label{si}
\ee
By contrast  with the mentioned Heisenberg double, this one is composed of two copies of RE algebras.

Also, we have introduced the matrix  $D=M^{-1}\hL$, whose entries are interpreted as $q$-partial derivatives in the entries of the matrix $M$. The matrix $D$ was shown to
satisfy the matrix relation \cite{GPS2}:
$$
R^{-1}D_1R^{-1}D_1=D_1 R^{-1} D_1R^{-1}.
$$
The permutation relations between the matrices $D$ and $M$ read as follows
\be
D_1RM_1R=RM_1R^{-1}D_1+ R\,1_B1_A.
\label{raz}
\ee
Below, we shall omit the unit factors $1_A$ and $1_B$ in order to simplify the formulas.

Note that our treatment of the entries of the matrix $D$ as analogs of partial derivatives is motivated by the following reason. If $R=P$, the relation (\ref{raz}) means that the entries
of the matrix $D$ are just  the usual partial derivatives in the  entries of the  matrix $M$, which commute with each other. In the general case, we call  the entries of the matrix $D$
{\em the quantum partial derivatives}.

Observe that a particular case of such a QD, corresponding to the Quantum Group $U_q(sl_2)$, was considered earlier in some papers (see  the references in \cite{GPS3}).

Let us go back to the QD $(\hLL(R), \MM(R),\sigma)$, where $\sigma$ is defined in (\ref{si}). It is convenient to make a linear shift of generators in the modified RE algebra
$\hLL(R)$ $L = I-(q-q^{-1})\hL$ and pass to the QD $(\LL(R), \MM(R),\sigma)$, where the permutation map for new generators takes the form:
$$
\sigma: R  L_1  R  M_1\,\to\, M_1  R  L_1  R^{-1}.
$$
The full list of the defining relations of this QD is as follows
$$
R  L_1  R  L_1=L_1  R  L_1  R, \qquad R  M_1  R  M_1=M_1  R  M_1  R,
$$
\be
R  L_1  R  M_1=M_1  R  L_1  R^{-1}.
\label{perm}
\ee
The equality (\ref{perm}) is a permutation relation between subalgebras $\LL(R)$ and $\MM(R)$.

This QD will play the main role in our considerations below. In particular we perform the spectral analysis of the Casimir operators of the subalgebra $\LL(R)$ acting
on the subalgebra $\MM(R)$ in accordance with the rule (\ref{AB-act}).

The technical tools for describing a representation of $\LL(R)$ in $\MM(R)$ are given by the so-called $R$-matrix representations of the Hecke algebras $H_k(q)$, $k\ge 1$.
Recall that a Hecke algebra $H_k(q)$ is generated by  elements $1, \tau_1,\dots, \tau_{k-1}$ subject to the relations
$$
\begin{array}{lcl}
\tau_i\,\tau_{i+1}\, \tau_i=\tau_{i+1}\,\tau_i\,\tau_{i+1},&\quad&1\le i\le k-2\\
\rule{0pt}{5mm}
\tau_i\, \tau_j=\tau_j\,\tau_i,&&|i-j|\geq 2\\
\rule{0pt}{5mm}
(q\, 1-\tau_i)(\qq\, 1+\tau_i)=0,&&q\in \C\setminus \{\pm 1, 0\}.
\end{array}
$$
A maximal commutative subalgebra in the Hecke algebra $H_k(q)$ is generated by the so-called {\it Jucys-Murphy} elements $J_i$, $1\le i\le k$, which are defined
by the recurrence relation (see \cite{OP} for detail):
$$
J_1 = 1,\quad J_i = \tau_{i-1}J_{i-1}\tau_{i-1}, \quad 2\le i\le k.
$$

With any Hecke symmetry $R$ we can associate the $R$-matrix representation of the Hecke algebra $H_k(q)$ in the space $V^{\otimes k}$ by the following
map
$$
\tau_i\,\mapsto\, R_i := I^{\ot (i-1)}\ot R\ot I^{\ot (k-i-1)}\in \mathrm{End}(V^{\otimes k}),\quad 1\le i \le k-1.
$$

Besides, we introduce the ``copies'' of the generating matrix $L$  of the RE algebra $\LL(R)$:
$$
L_{\ov 1}=L_{\un 1}=L_1,\qquad L_{\ov {i+1}}=R_i L_{\ov i}R_i^{-1},\qquad L_{\un {i+1}}=R_i^{-1}L_{\un i}R_i,\quad i\ge 1.
$$
The same notation will be used for the generating matrix $M$ of the RE algebra $\MM(R)$.

With this notation we can rewrite the permutation relation (\ref{perm}) using the $R$-matrix representation of the Jucys-Murphy elements:
$$
L_{\un 2} M_1= R_1^{-2} M_1 L_{\ov 2}=J_2^{-1} M_1L_{\ov 2}.
$$
This relation gives rise to the following  generalizations valid for a chain of the ``copies'' of the matrix $M$:
\be
L_{\un{k+1}} M_1M_{\ov 2}\dots M_{\ov k}=J_{k+1}^{-1} M_1M_{\ov 2}\dots M_{\ov k}L_{\ov{k+1}}.
\label{matr}
\ee

Now, introduce the counit $\varepsilon: \LL(R)\rightarrow {\Bbb C}$ in the algebra $\LL(R)$:
$$
\varepsilon (l_i^j)=\de_i^j, \quad \varepsilon (1 )=1_\C,\quad  \varepsilon (a\cdot b)=\varepsilon (a)\, \varepsilon (b),\quad \forall \, a,b \in \LL(R).
$$
According to the scheme described above we get the action of the algebra $\LL(R)$ onto $\MM(R)$:
\be
L_{\un{k+1}}\triangleright M_1\dots M_{\ov k}=J_{k+1}^{-1} M_1\dots M_{\ov k}.
\label{dva}
\ee
Hereafter, the notation $\triangleright $ means that in any entry of the matrix $L_{\un{k+1}} M_1\dots M_{\ov k}$ the element of the algebra $\LL(R)$ is applied as an
operator to the element of $\MM(R)$.

\begin{remark} \rm
In \cite{GPS2} the QD $(\LL(R),T(V))$ was considered. Here $T(V)$ is the free tensor algebra of the basic space $V$. The permutation relations
(see formula (3.20) in \cite{GPS2})  give rise  to the following action of generators of the algebra $\LL(R)$ onto the basis elements of $T(V)$:
\be
L_{\un{k+1}}\triangleright x_1\ot\dots \ot x_{ k}=J_{k+1}^{-1}x_1\ot\dots \ot x_{ k}.
\label{tri}
\ee
Note that this action is similar to (\ref{dva}). It is worth mentioning that representation (\ref{dva}) is defined by a {\it one-sided left} action of $L$ on a matrix $M_1\dots M_{\ov k}$
and due to this fact it is isomorphic to the direct sum of representations (\ref{tri}) in the spaces $V^{\ot k}$.
\end{remark}

In  \cite{IP} there was exhibited a way to associate a central element of the RE algebra  with any element of the Hecke algebra $H_k(q)$ (and even more, of the group algebra of
the braid group). We introduce some of these central elements.

First of all, we note that for a generic $q$ the Hecke algebra $H_k(q)$ is isomorphic to the group algebra $\C[S_k]$ of the symmetric group, therefore its representation
theory is the same as that of the symmetric group. In particular, to any partition $\la=(\la_1\geq \la_2\geq \dots \geq \la_k)$ of the integer $k$ we put into correspondence a
set of primitive idempotents $P_{(\la,a)}\in H_k(q)$ enumerated by the standard Young tables $(\lambda,a)$ constructed for the Young diagram of the partition $\lambda$.
Here the index $a$ indicates different standard tables $(\lambda ,a)$ with respect to some ordering. In the $R$-matrix representations these idempotents turn into orthogonal
projection operators which we denote as $P_{(\la, a)}(R)$.

According to \cite{IP} the central element of the algebra $\MM(R)$ corresponding to the projector $P_{(\lambda,a)}$ is given by the formula:
\be
s_{\la}(M)=\Tr_{R(1\dots k)}\,(P_{(\la, a)}(R)M_1M_{\ov 2}\dots M_{\ov k}).
\label{shur}
\ee
Hereafter, we use the following notation: $\Tr_{R(1\dots k)}=\Tr_{R(1)}\dots\Tr_{R(k)}$.

As was shown in \cite{OP}, the primitive idempotents $P_{(\la, a)}$ and $P_{(\la, a')}$, corresponding to different standard tables $(\lambda,a)$ and $(\lambda, a')$ of the same
Young diagram $\lambda$, are conjugated in the Hecke algebra. This means that there exists an invertible element $c(a,a')\in H_k(q)$ such that
$$
P_{(\la, a)}\, c(a,a')=c(a,a')\,P_{(\la, a')}.
$$
This property entails that the elements $s_{\la}(M)$ in (\ref{shur}) depend on the {\it diagrams} $\lambda$ but not on the on the tables. For this reason we omit the index $a$ in
their notation. We call the elements $s_{\la}(M)$ the $q$-Schur (or simply Schur) polynomials.

If $\la$ is a one-row or a one-column diagram, we get a {\em complete symmetric polynomial} $h_k(M)$ or {\em elementary symmetric polynomial} $e_k(M)$ respectively.
The corresponding projectors $S^{(k)}(R)$ and $A^{(k)}(R)$ are called the $R$-symmetrizer and the $R$-skew-symmetrizer.

Thus, the elementary symmetric polynomials are given by the formula
$$
e_k(M)=\Tr_{R(1\dots k)} (A^{(k)} M_1 M_{\ov 2}\dots M_{\ov k}),\quad k\ge 1.
$$
The set of $R$-skew-symmetrizers can be constructed with the use of the following recursion:
$$
A^{(1)}=I,\quad  A^{(k)}=\frac{1}{k_q}A^{(k-1)}\left(q^{k-1} I-(k-1)_q \, R_{k-1}\right)A^{(k-1)},\quad k\ge 2,
$$
where $ k_q=(q^k-q^{-k})/(q-\qq)$ is a $q$-integer. Note that $k_q\not=0$ $\forall k\in{\Bbb Z}$ for generic values of $q$.

Also, we need the central elements
$$
p_k(M)=\Tr_{R(1\dots k)}(R_{k-1}\dots R_2 R_1 \,M_1 M_{\ov 2}\dots M_{\ov k}\,)  =\Tr_R( M^k),\quad k\ge 1,
$$
which are called the {\it power sums}.

\begin{definition}\rm
The operators arising from all central elements of the algebras $\LL(R)$ or $\hLL(R)$ and acting in the spaces $V^{\ot k}$ in accordance with (\ref{tri}) or in the algebra $\MM(R)$
in accordance with (\ref{dva}) are called the $q$-Casimir operators.
\end{definition}

Note, that dealing with the $q$-Casimir operators of the RE algebra $\LL(R)$ is simpler in comparing with $q$-Casimir operators of the modified RE algebra
$\hLL(R)$. Recall, that any result obtained for the operators of $\LL(R)$ can be rewritten for the operators of $\hLL(R)$ by means of the linear shift of generators (\ref{iso}).

\section{Spectral analysis of the first Casimir operator}

In this section we perform the spectral analysis of the first $q$-Casimir and $q$-cut-and-join\footnote{Below we omit $q$ in the notation of these operators.} operators.
Note that the first cut-and-join operator and the first Casimir one coming from the algebra $\hLL(R)$ are equal to each other.

Let us intruduce the short-hand notation $\nu=q-\qq$. Besides, we shall systematically use the following property of the $R$-trace:
\be
\Tr_{R(k)} \,X_{\ov k}=\Tr_{R(k)} \,X_{\un k}= I_{1\dots k-1} \Tr_{R} \,X,
\label{syst}
\ee
where $X$ is an arbitrary $N\times N$ matrix and $I_{1\dots k-1}$ is the matrix of the identity operator in the space $V^{\otimes (k-1)}$.

\begin{proposition}
\label{prop:action}
The following relation takes place
$$
\Trr L\triangleright   M_1 M_{\ov 2}\dots M_{\ov k}=\left(q^{-N}N_q-\frac{\nu}{q^{2N}} \sum_{i=1}^k\, J_i^{-1}\right)\, M_1 M_{\ov 2}\dots M_{\ov k}.
$$
\end{proposition}
\medskip

\noindent
{\bf Proof.} By using (\ref{syst}), we get from  formula (\ref{dva}):
$$
\Trr L\triangleright  M_1 M_{\ov 2}\dots M_{\ov k}=\Tr_{R(k+1)}\,( J_{k+1}^{-1}\,  M_1M_{\ov 2}\dots M_{\ov k}).
$$
The Hecke condition (\ref{H-cond}) implies $R^{-1} = R-\nu I$, so we transform the element $J_{k+1}^{-1}$ to the sum of two terms:
$$
J_{k+1}^{-1} = R_{k}^{-1}J_k^{-1}R_k^{-1} = R_{k}^{-1}J_k^{-1}R_k - \nu R_{k}^{-1}J_k^{-1}.
$$
Then we apply (\ref{syst}) to the first term and take into account that\footnote{\label{foot-gen-rank} Recall that we are dealing with a $GL_N$ type Hecke symmetry $R$. For a
$GL_m$ type symmetry of an arbitrary rank $m$ the changes in formulas are minor: one should everywhere substitute $m$ instead of $N$.} (see \cite{GPS1})
$$
\Tr_{R(k+1)}R_k^{-1} = q^{-2N}I_{1\dots k}
$$
for the second term. So, we get:
$$
\Tr_{R(k+1)}\,( J_{k+1}^{-1}\,  M_1M_{\ov 2}\dots M_{\ov k}) = (\Tr_{R(k)}J_k^{-1} - \frac{\nu}{q^{2N}} J_k^{-1})\,M_1M_{\ov 2}\dots M_{\ov k}.
$$
Applying step by step the same transformations to $J_k^{-1}$, $J_{k-1}^{-1}$ and so on, we arrive at the following expression:
$$
\Tr_{R(k+1)}\,( J_{k+1}^{-1}\,  M_1M_{\ov 2}\dots M_{\ov k}) = \left(\Tr_{R(2)} J_2^{-1}-\frac{\nu}{q^{2N}} \sum_{i=2}^k\, J_i^{-1}\right) M_1 M_{\ov 2}\dots M_{\ov k}.
$$
As a consequence of the Hecke condition we have:
$$
J_2^{-1} = R_1^{-2} = I_{12} - \nu R_1^{-1}.
$$
Since for a $GL_N$ type Hecke symmetry $\Tr_RI = q^{-N}N_q$ (see \cite{GPS1}), we get
$$
\Tr_{R(2)} J_2^{-1} = \left(q^{-N}N_q-\frac{\nu}{q^{2N}}\right)I \equiv  q^{-N}N_q I -\frac{\nu}{q^{2N}}J_1^{-1},
$$
which completes the proof.\hfill\rule{6.5pt}{6.5pt}

To formulate the next proposition we need the notion of the {\it content} of a box of a Young diagram. Take a Young diagram corresponding
to some partition $\lambda$ and consider an arbitrary box of the diagram located at the intersection of the $m$-th row and the $n$-th column.
The {\it content } $c$ of the box is an integer number $c = n-m$.

\begin{proposition} \label{four}
Let $\lambda\vdash k$ be a partition and $(\lambda,a)$ be a standard Young table corresponding to the Young diagram of the partition $\lambda$.
Then the following relation holds
$$
\Tr_R L\triangleright P_{(\la,a)}(R) M_1M_{\ov 2}\dots M_{\ov k} =\chi_{\la} (\Tr_R L)P_{(\la,a)}(R)M_1M_{\ov 2}\dots M_{\ov k},
$$
where\rm
\be
\chi_{\la} (\Trr L)=\frac{N_q}{q^{N}}-\frac{\nu}{q^{2N}} \sum_{i=1}^k  \, q^{-2c(i)},
\label{dec}
\ee\it
the sum is taken over all boxes of the table $(\la, a)$ and $c(i)$ is the content of the box in which the integer $i$ is located.
\end{proposition}

Note, that in (\ref{dec}) the sum is taken over {\it all} boxes of  $(\la,a)$ and therefore it depends only on the {\it diagram} $\lambda$ but not
on the table $(\lambda, a)$. For this reason we omit the index $a$ in the notation of the character $\chi_{\la}$.
\medskip

\noindent
{\bf Proof.} Take into account that due to the braid relarion on $R$ we have
$$
\begin{array}{l}
R_1^{\pm 1}R_2^{\pm 1}\dots R_k^{\pm 1} R_{i} = R_{i+1} R_1^{\pm 1}R_2^{\pm 1}\dots R_k^{\pm 1} \qquad  1\le \forall\,i\le k-1,\\
\rule{0pt}{6mm}
R_k^{\pm 1}\dots R_2^{\pm 1}R_1^{\pm 1} R_j = R_{j-1} R_k^{\pm 1}\dots R_2^{\pm 1}R_1^{\pm 1}\qquad 2\le \forall \, j\le k.
\end{array}
$$
Therefore, the matrix $L_{\underline{k+1}} = R_{k}^{-1}\dots R_1^{-1}L_1R_1\dots R_k$ commute with any polynomial $Q$ in the matrices $R_i$, $1\le i\le k-1$:
$$
L_{\underline{k+1}} Q(R_1,R_2, \dots, R_{k-1}) = Q(R_1,R_2, \dots, R_{k-1}) L_{\underline{k+1}}.
$$
In the $R$-matrix representation any primitive idempotent $P_{(\la, a)}\in H_k(q)$ is given by a projector $P_{(\lambda,a)}(R)$ which is a polynomial in the matrices $R_i$,
$1\le i\le k-1$ and, consequently, commute with $L_{\underline{k+1}}$. So, taking into account (\ref{dva}) we find:
\be
L_{\un{k+1}}\triangleright  P_{(\la, a)}(R)  M_1M_{\ov 2}\dots M_{\ov k}=
P_{(\la, a)}(R) J_{k+1}^{-1} M_1M_{\ov 2}\dots  M_{\ov k}.
\label{action-L}
\ee
Besides, it can be shown (see \cite{OP}) that
$$
J_i P_{(\la,a)}(R)=P_{(\la,a)}(R)J_i=q^{2c(i)}P_{(\la,a)}(R).
$$
So, on applying the $R$-trace to the left and right hand sides of (\ref{action-L}) and on taking into account Proposition \ref{prop:action} we come to the desired result
(\ref{dec}).\hfill \rule{6.5pt}{6.5pt}

\medskip

Let us observe that Proposition \ref{four} remains true if we replace the product $M_1M_{\ov 2}\dots M_{\ov k}$ by basis vectors $x_1\otimes \dots \otimes x_k$ of
the space $V^{\ot k}$ (see (\ref{tri})).

\begin{corollary}
The Schur polynomials $s_{\la}(M)$ {\rm (\ref{shur})} are eigenvectors of the Casimir operator $\Tr_R L$ with respect to the action {\rm (\ref{dva})}:
$$
\Tr_R L\triangleright s_\lambda(M) = \chi_\lambda(\Tr_RL)\,s_\lambda(M),
$$
where the eigenvalue $\chi_{\la} (\Tr_R L)$ is defined in {\rm (\ref{dec})}.
\end{corollary}

\begin{definition} \rm If the Schur polynomial $s_{\la}(M)$ is an eigenvector of a given Casimir operator $\CC$
$$
\CC\triangleright s_{\la}(M)=\chi_{\la}(\CC)\, s_{\la}(M),
$$
then the corresponding eigenvalue $\chi_{\la}(\CC)$ is called  the $\la$-character of $\CC$.
\end{definition}

\begin{remark}\rm
Note that for a $GL_N$ type Hecke symmetry $R$ the projectors $P_{(\lambda,a)}(R)$ are identically equal to zero if the partition $\lambda$
contains more than $N$ nonzero components. So in what follows we will consider only partitions $\lambda$ with $\lambda_{i}=0$ for $i\ge N+1$
(see also footnote \ref{foot-gen-rank}).
\end{remark}

\begin{proposition}\label{prop:8}
Let $\lambda = (\la_1,\la_2,\dots,\la_N)$ be a partition. Then the $\la$-character {\rm (\ref{dec})} of the Casimir element $\Tr_R L$ can be written in the following form
\rm
\be
\chi_{\la}(\Trr L)=\qq\,\sum_{k=1}^N  q^{-2(\la_k+N-k)}.
\label{char1-mu}
\ee
\end{proposition}
\medskip

\noindent
{\bf Proof.} Let us multiply formula  (\ref{dec}) by $q$ and calculate the sums of terms $q^{-2c(i)}$ corresponding to each row of the Young diagram $\lambda$.
For such a partial sum corresponding to a $j$-th row we have:
$$
\nu\, q^{1-2N} \left( q^{2(j-1)}+q^{2(j-2)}+\dots +q^{2(j-\la_j)}\right)=q^{-2(N-j)}-q^{-2(\la_j+N-j)}.
$$
So, the second term of  (\ref{dec}), multiplied by $q$,  is equal to the expression
\be 
\nu \, q^{1-2N}\, \sum_{i=1}^k q^{-2c(i)}=\sum_{j=1}^s\left(q^{-2(N-j)}-q^{-2(\la_j+N-j)}\right), \label{pro}
\ee
where $s$ is the number of non-trivial components of the partition $\la = (\la_1,\la_2,\dots,\la_s,0,\dots, 0)$, $\la_s>0$.

Taking into account that
$$
q^{1-N}N_q=1+q^{-2}+\dots + q^{-2(N-1)},
$$
we arrive at the final result:
$$
q\,\chi_\lambda(\Tr_RL) = 1+q^{-2}+\dots +q^{-2(N-s-1)}+\sum_{j=1}^s q^{-2(\la_j+N-j)} = \sum_{k=1}^N  q^{-2(\la_k+N-k)},
$$
since $\lambda_k = 0$ for $s+1\le k\le N$.\hfill \rule{6.5pt}{6.5pt}
\medskip

Now, we introduce the first quantum cut-and-join operator setting by definition:
$$
W^{(1)} =\Trr \hL = \Tr_R MD,
$$
where the matrix $D$, composed of the  quantum partial derivatives, was defined in section 2 (see, in particular, (\ref{raz})). Since the entries of the matrix $\hL= MD$ generate the
modified RE algebra $\hLL(R)$, we are able to get the $\la$-characters of the operator $W^{(1)}=\Tr_R \hL$ by using formulae (\ref{iso}) and (\ref{dec}).
\begin{proposition}
For any partition $\lambda=(\lambda_1,\lambda_2,\dots,\lambda_N)$, $\lambda\vdash k\in {\Bbb Z}_+$, the following relation holds:
$$
W^{(1)}\triangleright s_{\la}(M) =\chi_{\la}(\Tr_R \hL) s_{\la}(M),
$$
where the $\lambda$-character of the operator $W^{(1)}$ reads
$$
\chi_{\la} (\Trr \hL)=\frac{1}{q^{2N}} \sum_{i=1}^k  \, q^{-2c(i)} = \frac{1}{q^{2N}} \sum_{i=1}^N
q^{-(\lambda_i +1-2i)}(\lambda_i)_q.
$$
\end{proposition}

{\bf Proof.} The first equality follows from (\ref{iso}), (\ref{dec}), and the relation $\Trr I=q^{-N}\, N_q$. The second equality follows from (\ref{pro}).  \hfill \rule{6.5pt}{6.5pt}

\section{Higher Casimir and cut-and-join operators}

In this section we shall consider the general case of a $GL_m$ type Hecke symmetry with an arbitrary rank $m\le N$. For such a symmetry $R$ the generating matrix $L$ of
the algebra $\LL(R)$ is subject to the matrix Cayley-Hamilton (CH) identity (see for instance \cite{GPS1} and the references therein):
\be
L^m-q e_1(L) L^{m-1}+ q^2e_2(L) L^{m-2}+\dots  +(-q)^{m-1} e_{m-1}(L) L+(-q)^{m} e_{m}(L)\, I=0.
\label{CH}
\ee
The corresponding {\it characteristic polynomial} $Q(t)$ meeting the relation $Q(L)=0$ is of the form:
$$
Q(t)= \sum_{k=0}^m\, t^{m-k}(-q)^ke_k(L).
$$
Its roots $\{\mu_i\}_{1\le i\le m}$ are called the {\em eigenvalues} of the matrix $L$. These roots are elements of a central extension of the algebra $\LL(R)$.
Thus, we have
\be
 q^k e_k(L)=\sum_{1\leq i_1<\dots <i_k\leq m} \mu_{i_1}\dots \mu_{i_k}.
\label{elem-mu}
\ee
In particular, $\Tr_R L=e_1(L)=\qq\, \sum_{k=1}^m\, \mu_k$.

It is tempting to assign to the eigenvalues $\mu_i$ their $\la$-characters $\chi_\la(\mu_i)$ whose values would be consistent with the relations (\ref{elem-mu}).

\begin{conjecture} \label{conj:10}
Given a partition $\lambda = (\lambda_1,\lambda_2,\dots,\lambda_m)$, we set
$$
\chi_\la(\mu_i)=q^{-2(\la_i+m-i)}.
$$
This assignement is consistent with the relations {\rm (\ref{elem-mu})}, that is
$$
\chi_{\la}(e_k(L))=q^{-k}\,\sum_{1\leq i_{1}<\dots < i_k \leq m}   \chi_{\la}(\mu_{i_1})\dots \chi_{\la}(\mu_{i_k}).
$$
\end{conjecture}

In Proposition \ref{prop:8}  this conjecture was proved for the particular case $k=1$. Also, we checked this conjecture in a few low-dimensional cases.

\begin{remark} \rm As was shown by straightforward calculations, the $\la$-character of the Casimir operator $\Trr L$ does not depend on the Young table $(\lambda,a)$. So far
we do not have analogous explicit calculations for the characters of the higher Casimir operators. Nevertheless, it is not difficult to see that the $\la$-characters of the these operators
do not depend on the table indeed. It is due to the property, mentioned above, that the idempotents $P_{(\la,a)}$, corresponding to different tables, are conjugated to each other.
\end{remark}

Assuming Conjecture \ref{conj:10} to be true, we become able to compute the $\la$-characters of the elementary symmetric polynomials and consequently of the power
sums since these two families are related by $q$-analogs of the Newton relations. The reader is referred to \cite{GPS4} for more details on the relations between these families.

Note that the matrix $\hL$ is also subject to a version of the CH identity. It can be readily obtained from (\ref{CH}) and (\ref{iso}). Let $\hmu_i$ be the eigenvalues of the characteristic polynomial for the matrix $\hL$. They are related to the eigenvalues of the matrix $L$ by the formula
$$
\mu_i=1-\nu \hmu_i.
$$
Thus, we get
$$
\chi_\la({\hat{\mu}_k})=\frac{1- q^{-2\, (\la_k+m-k)}}{q-\qq} = q^{-(\lambda_k+m-k)}(\lambda_k+m-k)_q,\quad 1\le k \le m.
$$
If $m=N$ and $R\to P$ as $q\to 1$, then by passing to the limit $q\rightarrow 1$ in the above relation we obtain the $\la$-characters of the eigenvalues $\hat{\mu}_k$ of the matrix $\hL$ generating the algebra $U(gl(N))$:
$$
\chi_\la({\hat{\mu}_k})=\la_k+N-k,\quad 1\le k\le N.
$$
Note that a similar result is also valid for the generating matrix $\hL$ of the modified RE algebra, associated with any involutive symmetry $R$
(i.e. $R^2=I$), provided $R$ is the $q\to 1$ limit of a $GL_N$ type Hecke symmetry.

Now, we pass to construction of the higher cut-and-join operators $W^\De$, which are similar to those  (\ref{cut}) but with a new meaning of the matrix $\hL$, the traces and
the normal ordering. In our setting the matrix $\hL$, giving rise to the cut-and-join operators, equals to the matrix product $MD$, where $D$ and $M$ are the generating matrices
of the algebras $\LL(R^{-1})$ and $\MM(R)$ respectively, their permutation relations are given by (\ref{raz}).

Thus, in order to introduce the quantum analogs of the operators $W^{\De}$  we have to define the quantum normal ordering rule. For this purpose we replace each matrix $\hL$
in the definition of $W^{\De}$ by the product $MD$ and then, in full analogy with the classical pattern, we push all the matrices $D$ through all the matrices $M$ to the most right
position by means of the following permutation relations:
\be
D_1RM_1R=RM_1 R^{-1}D_1\quad \Leftrightarrow \quad D_1RM_1=RM_1R^{-1}D_1R^{-1}.
\label{razz}
\ee
These permutation relations are obtained from (\ref{raz}) by omitting the last (constant) term. Thus, we define  the $q$-analog of the normal ordering for
the product of $q$-differential operators.

Besides, similarly to the classical case considered in \cite{MMZ}, the cut-and-join operators thus defined can be expressed in terms of the quantum Casimir ones.

 Let us consider a few examples. Namely, we compute the operators $W^\De$ for $\De=(2,0,\dots,0)$, $\De=(1,1,0,\dots,0)$, and $\De=(3,0,\dots,0)$.

\begin{example} \rm
For $\Delta = (2,0,\dots,0)$ the cut-and-join operator reads $W^{(2)} = \,:\!\!\Tr_R\hL^2\!\!:$. In order to get the explicit expression for the normal
ordered form, we use the identity
$$
\Tr_R\hL^2 = \Tr_{R(12)}(\hL_1\hL_{\ov 2} R_1)=\Tr_{R(12)}(\hL_1R_1\hL_1) = \Tr_{R(12)}(M_1D_1R_1M_1D_1)
$$
and apply the permutation relations (\ref{razz}):
\begin{eqnarray*}
W^{(2)} = :\Trr \hL^2:\!\!\!&=&\!\!\!:\Tr_{R(12)} (M_1{\un{D_1R_1M_1}}D_1):=\Tr_{R(12)}(M_1R_1M_1R_1^{-1}D_1R_1^{-1}D_1)\\
&=& \Tr_{R(12)}(M_1M_{\ov 2}\,D_1R_1^{-1}D_1)=\Tr_{R(12)}(M_1M_{\ov 2}D_{\ov 2} D_1R_1^{-1}).
\end{eqnarray*}
A straightforward computation leads to the following connection between $W^{(2)}$ and the Casimir operators:
\be
W^{(2)} =\Tr_R \hL^2-q^{-m}m_q\Tr_R\hL.
\label{mor1}
\ee
\end{example}

\begin{example}\rm
For $\Delta = (1,1,0,\dots,0)$ we have the cut-and-join operator
$$
W^{(1,1)} = : (\Tr_R\hL)^2: = :\Tr_{R(12)}(\hL_1\hL_{\ov 2}):.
$$
Its explicit form is
\begin{eqnarray*}
W^{(1,1)}\!\!\!&=&\!\!\!:\Tr_{R(12)} (M_1{\un{D_1R_1M_1}}D_1R_1^{-1}): = \Tr_{R(12)}(M_1R_1M_1R_1^{-1}D_1R_1^{-1} D_1R_1^{-1} )\\
&=&\Tr_{R(12)}(M_1M_{\ov 2}D_{\ov 2}D_1R_1^{-2}).
\end{eqnarray*}
 Also, we get the following expression of the operator $W^{(1,1)}$ in terms of the Casimir operators:
 \be
W^{(1,1)}=(\Tr_R \hL)^2-q^{-2m} \,\Tr_R \hL.
\label{mor2}
\ee
\end{example}

\begin{example}\rm
At last, we give the explicit form of the cut-and-join operator $W^{(3)} =:\Tr_R\hL^3:$, corresponding to $\De=(3,0,\dots,0)$. To this end we use the identity
$$
\Tr_R\hL^3 = \Tr_{R(123)}(\hL_1\hL_{\ov 2}\hL_{\ov 3}R_2R_1)
$$
and then apply the relations (\ref{razz}). Omitting the straightforward calculations we give the answer
$$
W^{(3)} = :\Tr_R(\hL^3): = \Tr_{R(123)}(M_1M_{\ov 2}M_{\ov 3}D_{\ov 3}D_{\ov 2}D_1(R_2^{-1} R_1^{-1})^2).
$$
The expression of $W^{(3)}$ via the Casimir operators is as follows
\be
W^{(3)} =\Tr_R \hL^3-2\frac{m_q}{q^m} \,\Tr_R \hL^2-\frac{1}{q^{2m}}\, (\Tr_R \hL)^2+(q^{-4m}+q^{-2m}m_q^2\,)\, \Tr_R \hL.
\label{mor3}
\ee
\end{example}

At last, we point out that the $\la$-characters of all operators under consideration do not depend on the concrete  form   of the initial Hecke symmetry $R$, but only on its rank.

As we said above, the classical cut-and-join operators are usually expressed in terms of the power sums $p_k = \Tr(M^k)$, $k\ge 1$, and the partial derivatives in $p_k$. In our
setting such an expression is impossible and the action of the cut-and-join operators on the power sums $p_k$ (and functions in $p_k$) has to be computed directly. We plan to consider this problem in our subsequent publications.


\begin{thebibliography}{GPS3}

\bibitem[Go]{Go} Goulden I.P., A differential operator for symmetric functions and the combinatorics of multiplying transpositions, Transactions of AMS 344 (1994) 421--440.

\bibitem[G]{G} Gurevich D.,  Algebraic aspects of the Yang-Baxter equation,  Leningrad Math. J.  2 (1991) 801 -- 828.


\bibitem[GPS1]{GPS1} Gurevich D., Pyatov P., Saponov P.,  Representation theory of (modified) reflection equation algebra of $GL(m|n)$ type,
 St. Petersburg Math. J. 20 (2008) 213 -- 253.

\bibitem[GPS2]{GPS2} Gurevich D., Pyatov P., Saponov P., Braided Differential operators on quantum algebras, J. of Geometry and Physics 61 (2011) 1485 -- 1501.

\bibitem[GPS3]{GPS3} Gurevich D., Pyatov P., Saponov P., Braided Weyl algebras and differential calculus on $U(u(2))$, J. of Geometry and
 Physics 62 (2012) 1175 -- 1188.

\bibitem[GPS4]{GPS4} Gurevich D., Pyatov P., Saponov P.,  Spectral parametrization for power sums of quantum supermatrices, Theor. Math. Phys. 159 (2009) 587 -- 597.

\bibitem[GS1]{GS1} Gurevich D., Saponov P., From Reflection Equation Algebra to Braided Yangians, Proceedings of the 1st International Conference on Mathematical Physics,
Grozny, Russia, 2016. Springer Proceedings in Math. and Statistics V.273 (2018).

\bibitem[GS2]{GS2} Gurevich D., Saponov P.,  Doubles of Associative algebras and their Applications, Phys. of Particles and Nuclei Letters 17, no. 5 (2020) 774 -- 778.

\bibitem[IP]{IP} Isaev A., Pyatov P., Spectral Extension of the Quantum Group Cotangent Bundle, Comm. Math. Phys. 288 (2009) 1137 -- 1179.

\bibitem[MMN]{MMN} Mironov A., Morosov A., Natanzon S., Complete set of Cut-and-Join operators in Hurwitz-Kontsevich theory, Theor. Math. Phys. 166 (2011) 1 -- 22.

\bibitem[MMZ]{MMZ} Mironov A., Morosov A., Zhabin A., Connection between cut-and join and Casimir operators, ArXiv: 2105.10978.

\bibitem[OP]{OP} Ogievetsky O., Pyatov P., Lecture on Hecke algebras  In “Symmetries and Integrable systems”, Dubna publishing, 2000,
Preprint CPT-2000/P.4076.

\bibitem[PP]{PP} Perelomov A., Popov V., Casimir operators for the classical groups, Dokl. Akad. Nauk SSSR, 1967, 174, 287--290

\bibitem[S]{S} Saponov P., Weyl approach to representation theory of reflection equation, Journal of Physics A: Mathematical and General,
37, no. 18 (2004) 5021 -- 5046.

\end{thebibliography}
\end{document}